# PRELIMINARIES ON PSEUDONTRACTIONS IN THE INTERMEDIATE SENSE FOR NON-CYCLIC AND CYCLIC SELF-MAPPINGS IN METRIC SPACES

M. De la Sen

Instituto de Investigacion y Desarrollo de Procesos, Universidad del Pais Vasco

Campus of Leioa (Bizkaia) – Aptdo. 644- Bilbao, 48080- Bilbao, SPAIN, email: *manuel.delasen@ehu.es*

**Abstract**. A contractive condition is addressed for extended 2-cyclic self-mappings on the union of a finite number of subsets of a metric space which are allowed to have a finite number of successive images in the same subsets of its domain. It is proven that if the space is uniformly convex and the subsets are non-empty, closed and convex then all the iterations converge to a unique closed limiting finite sequence which contains the best proximity points of adjacent subsets and reduce to a unique fixed point if all such subsets intersect.

## 1. Introduction

Strict pseudocontractive mappings and pseudocontractive mappings in the intermediate sense formulated in the framework of Hilbert spaces have received a certain attention in the last years concerning their convergence properties and the existence of fixed points (see, for instance, [1-4] and references therein). Results about the existence of a fixed point are discussed in those papers. On the other hand, important attention has been paid during the last decades to the study of the convergence properties of distances in cyclic contractive self-mappings on $p$ subsets $A_i \subset X$ of a metric space $(X,d)$, or a Banach space $(X,\|\ \|)$. The cyclic self-mappings under study have been of standard contractive or weakly contractive types and of Meir-Keeler type. The convergence of sequences to fixed points and best proximity points of the involved sets has been investigated in the last years. See, for instance, [6-21] and references therein. This research investigates the convergence properties and the existence of fixed points of a generalized version of pseudocontractve, strict pseudocontractive and asymptotically pseudocontractive in the intermediate sense in the more general framework of metric spaces. The case of 2-cyclic pseudocontractive self-mappings is also considered. The combination of constants defined the contraction may be different on each of the subsets and only the product of all the constants is requested to be less than unity. On the other hand, the self-mapping can perform a number of iterations on each of the subsets before transferring its image to the next adjacent subset of the 2-cyclic self-mapping. The existence of a unique closed finite limiting sequence on any sequence of iterations from any initial point in the union of the subsets is proven if $X$ is a uniformly convex Banach space and all the subsets of $X$ are nonempty, convex and closed. Such a limiting sequence is of size $q \geq p$ (with the inequality being strict if there is at least one iteration with image in the same subset as its domain) where $p$ of its elements (all of them if $q = p$) are best proximity points between adjacent subsets. In the case that all the subsets $A_i \subset X$ intersect, the above limit sequence reduces to a unique fixed point allocated within the intersection of all such subsets.

## 2. Contractions and pseudocontractions in the intermediate sense

If $H$ is a real Hilbert space with an inner product $\langle .,. \rangle$ and a norm $\|.\|$ and $A$ is a nonempty closed convex subset of $H$ then the self-mapping $T: A \to A$ is said to be an asymptotically $\beta$-strictly pseudocontractive self-mapping in the intermediate sense for some $\beta \in [0,1)$ if

$$\limsup_{n \to \infty} \sup_{x,y \in C} \left( \|T^n x - T^n y\|^2 - \alpha_n \|x-y\|^2 - \beta \|(I-T^n)x - (I-T^n)y\|^2 \right) \leq 0, \quad \forall x,y \in A \qquad (2.1)$$

for some sequence $\{\alpha_n\} \subset [1,\infty)$, $\alpha_n \to 1$ as $n \to \infty$, [1-5]. Such a concept was firstly introduced in [1]. If (2.1) holds for $\beta = 1$ then $T: A \to A$ is said to be an asymptotically pseudocontractive self-mapping in the intermediate sense. Finally, if $\alpha_n \to \alpha \in [0,1)$ as $n \to \infty$ then $T: A \to A$ is asymptotically $\beta$-strictly contractive in the intermediate sense, respectively, asymptotically contractive in the intermediate sense if $\beta = 1$. If (2.1) is changed to the stronger condition:

$$\left( \|T^n x - T^n y\|^2 - \alpha_n \|x-y\|^2 - \beta \|(I-T^n)x - (I-T^n)y\|^2 \right) \leq 0; \quad \forall x,y \in A, n \in N \qquad (2.2)$$

then the above concepts translate into $T: A \to A$ being an asymptotically $\beta$-strictly pseudocontractive self-mapping, an asymptotically pseudocontractive self-mapping and asymptotically contractive, respectively. Note that (2.1) is equivalent to:

$$\|T^n x - T^n y\|^2 \leq \alpha_n \|x-y\|^2 + \beta \|(I-T^n)x - (I-T^n)y\|^2 + \xi_n; \quad \forall x,y \in A, \forall n \in N \qquad (2.3)$$

or, equivalently,

$$\langle T^n x - T^n y, x-y \rangle \leq \frac{1}{2\beta} \left[ (\alpha_n + \beta)\|x-y\|^2 + (\beta-1)\|T^n x - T^n y\|^2 + \xi_n \right]; \quad \forall x,y \in A, n \in N \qquad (2.4)$$

where

$$\xi_n := \max \left\{ 0, \sup_{x,y \in C} \left( \|T^n x - T^n y\|^2 - \alpha_n \|x-y\|^2 - \beta \|(I-T^n)x - (I-T^n)y\|^2 \right) \right\}; \quad \forall n \in N \qquad (2.5)$$

Note that the high-right-hand-side term $\|(I-T^n)x - (I-T^n)y\|^2$ of (2.3) is expanded as follows for any $x,y \in C$:

$$\|x-y\|^2 + \|T^n x - T^n y\|^2 - 2\|x-y\|\|T^n x - T^n y\|$$

$$\leq \|(I-T^n)x - (I-T^n)y\|^2 \leq \langle x - T^n x, y - T^n y \rangle^2 = \langle x-y, T^n x - T^n y \rangle^2$$

$$= \|x-y\|^2 + \|T^n x - T^n y\|^2 + 2\langle T^n x - T^n y, x-y \rangle = \langle x-y, T^n x - T^n y \rangle \langle x-y, T^n x - T^n y \rangle$$

$$\leq \|x-y\|^2 + \|T^n x - T^n y\|^2 + 2|\langle T^n x - T^n y, x-y \rangle| \leq \|x-y\|^2 + \|T^n x - T^n y\|^2 + 2\|x-y\|\|T^n x - T^n y\|$$

(2.6)

The objective of this paper is to discuss the various pseudocontractive in the intermediate sense concepts in the framework of metric spaces endowed with an homogeneous and translation-invariant metric and also to generalize them to the $\beta$ - parameter to eventually be replaced with a sequence $\{\beta_n\}$ in $(0,1)$.

Now, if instead of a real Hilbert space $H$ endowed with an inner product $\langle .,. \rangle$ and a norm $\|.\|$, we deal with any generic Banach space $(X, \|.\|)$ then its norm induces an homogeneous and translation invariant metric $d: X \times X \to \mathbf{R}_{0+}$ defined by $d(x,y) = d(x-y, 0) = \|x-y\|^{1/2}$; $\forall x, y \in A$ so that (2.6) takes the form:

$$d^2(x,y) + d^2(T^n x, T^n y) - 2d(x,y)d(T^n x, T^n y)$$
$$\leq \|(I-T^n)x - (I-T^n)y\|^2 = d^2(x-y - (T^n x - T^n y), 0) = d^2(x-y, T^n x - T^n y)$$
$$\leq (d(x-y, 0) + d(T^n x - T^n y, 0))^2 = (d(x,y) + d(T^n x, T^n y))^2$$
$$= d^2(x,y) + d^2(T^n x, T^n y) + 2d(x,y)d(T^n x, T^n y); \quad \forall x, y \in A \qquad (2.7)$$

Define

$$\mu_n(x,y) := \min\left(\rho \in [-1,1]: d^2(x-y, T^n x - T^n y) \leq d^2(x,y) + d^2(T^n x, T^n y) + 2\rho d(x,y)d(T^n x, T^n y)\right)$$
$$; \forall x, y \in A, \forall n \in \mathbf{N} \qquad (2.8)$$

which exists since it follows from (2.7), since the metric is homogeneous and translation-invariant, that:

$$\{1\} \subset \left\{\rho \in \mathbf{R}: \|(I-T^n)x - (I-T^n)y\|^2 \leq d^2(x,y) + d^2(T^n x, T^n y) + 2\rho d(x,y)d(T^n x, T^n y)\right\} (\neq \emptyset) \quad (2.9)$$

The following result holds related to the discussion (2.7)-(2.9) in metric spaces:

**Theorem 2.1**. Let $(X, d)$ be a metric space and consider a self-mapping $T: X \to X$. Assume that the constraint below holds:

$$d^2(T^n x, T^n y) \leq \alpha_n(x,y)d^2(x,y) + \beta_n(x,y)(d^2(x,y) + d^2(T^n x, T^n y))$$
$$+ 2\mu_n(x,y)\beta_n(x,y)d(x,y)d(T^n x, T^n y) + \xi_n(x,y); \quad \forall x, y \in X, \forall n \in \mathbf{N} \qquad (2.10)$$

with

$$\xi_n = \xi_n(x,y)$$
$$:= \max\left(0, (1-\beta_n(x,y))d^2(T^n x, T^n y) - (\alpha_n(x,y) + \beta_n(x,y))d^2(x,y) - 2\mu_n(x,y)\beta_n(x,y)d(x,y)d(T^n x, T^n y)\right)$$
$$\to 0 \quad ; \forall x, y \in X \text{ as } n \to \infty \qquad (2.11)$$

for some parameterizing bounded real sequences $\{\alpha_n(x,y)\}$, $\{\beta_n(x,y)\}$ and $\{\mu_n(x,y)\}$ of general terms $\alpha_n = \alpha_n(x,y)$, $\beta_n = \beta_n(x,y)$, $\mu_n = \mu_n(x,y)$ satisfying the following constraints:

$$\left[\left(\mu_n(x,y) \in \left[-\frac{\alpha_n(x,y) + \beta_n(x,y)}{2\beta_n(x,y)}, \frac{1-\alpha_n(x,y) - 2\beta_n(x,y)}{2\beta_n(x,y)}\right]\right) \wedge (\beta_n(x,y) < 1)\right]$$
$$\vee \left[\left(\mu_n(x,y) < -\frac{\alpha_n(x,y) + \beta_n(x,y)}{2\beta_n(x,y)}\right) \wedge (\beta_n(x,y) > 1) \Leftrightarrow \xi_n(x,y) = 0\right]$$
$$\vee \left[\mu_n(x,y) \in \left[\frac{1-\alpha_n(x,y) - 2\beta_n(x,y)}{2\beta_n(x,y)}, \frac{1-\beta_n(x,y)}{2\beta_n(x,y)}\right]\right] \quad ; \forall x, y \in X, \forall n \in \mathbf{N} \qquad (2.12)$$

with $\limsup\limits_{n\to\infty} [\beta_n(x,y)\max(1, 1+2\mu_n(x,y))] < 1$ and, furthermore, the following limit exists:

$$\left(\mu_n(x,y) - \frac{1-\alpha_n(x,y)-2\beta_n(x,y)}{2\beta_n(x,y)}\right) \to 0; \quad \forall x,y \in X \text{ as } n \to \infty$$

$$\Leftrightarrow \alpha_n + 2\beta_n(1+\mu_n) \to 1; \quad \forall x,y \in X \text{ as } n \to \infty \tag{2.13}$$

Then, the following properties hold:

**(i)** $\exists \lim\limits_{n\to\infty} d(T^n x, T^n y) \leq d(x,y)$ for any $x,y \in X$ so that $T: X \to X$ is asymptotically nonexpansive.

**(ii)** Let $(X,d)$ be complete, $d: X \times X \to \mathbf{R}_{0+}$ be, in addition, a translation-invariant homogeneous norm and $(X, \|\ \|) \equiv (X,d)$, with $\|\ \|$ being the metric-induced norm from $d: X \times X \to \mathbf{R}_{0+}$, be a uniformly convex Banach space. Assume also that $T: X \to X$ is continuous. Then, any sequence $\{T^n x\}$; $\forall x \in A$ is bounded and convergent to some point $z_x = z_x(x) \in C$, being in general dependent on $x$, in some nonempty bounded, closed and convex subset $C$ of $A$, where $A$ is any nonempty bounded subset of $X$. Also, $d(T^n x, T^{n+m} x)$ is bounded; $\forall n, m \in N$, $\lim\limits_{n\to\infty} d(T^n x, T^{n+m} x) = 0$; $\forall x \in A$, $\forall m \in N$ and $z_x = z_x(x) = Tz_x \in C$ is a fixed point of the restricted self-mapping $T: C \to C$; $\forall x \in A$. Furthermore,

$$\lim\limits_{n\to\infty} \left(d^2(T^{n+1}x, T^{n+1}y) - d^2(T^n x, T^n y)\right) = 0; \forall x, y \in A. \tag{2.14}$$

*Proof*: Consider two possibilities for the constraint (2.10), subject to (2.11), to hold for each given $x, y \in X$ and $n \in N$ as follows:

**A)** $d(T^n x, T^n y) \leq d(x,y)$ for any $x,y \in X$, $n \in N$. Then, one gets from (2.10):

$d^2(T^n x, T^n y) \leq (\alpha_n + \beta_n)d^2(x,y) + \beta_n d^2(T^n x, T^n y) + 2\mu_n \beta_n d^2(x,y) + \xi_n$

$$\Rightarrow d(T^n x, T^n y) \leq k_{an} d^2(x,y) + \frac{\xi_n}{1-\beta_n} \tag{2.15}$$

; $\forall x,y \in A$, $\forall n \in N$, where

$$k_{an} = k_{an}(x,y) = \frac{\alpha_n + \beta_n(1+2\mu_n)}{1-\beta_n} \to 1; \forall x,y \in X \text{ as } n \to \infty \tag{2.16}$$

which holds from (2.12)-(2.13) if $\limsup\limits_{n\to\infty} \beta_n(x,y) < 1$ since $\left(\mu_n(x,y) - \frac{1-\alpha_n(x,y)-2\beta_n(x,y)}{2\beta_n(x,y)}\right) \to 0$;

$\forall x,y \in X$ as $n \to \infty$ in (2.13) is equivalent to (2.16). Note that $0 \leq k_{an} \leq 1$ is ensured either with $\min(\alpha_n + \beta_n(1+2\mu_n), 1-\beta_n) \geq 0$ or with $\max(\alpha_n + \beta_n(1+2\mu_n), 1-\beta_n) \leq 0$ if

$$\left[\left(\mu_n(x,y) \in \left[-\frac{\alpha_n(x,y)+\beta_n(x,y)}{2\beta_n(x,y)}, \frac{1-\alpha_n(x,y)-2\beta_n(x,y)}{2\beta_n(x,y)}\right]\right) \wedge (\beta_n(x,y) \in (0,1))\right]$$

$$\vee \left[\left(\mu_n(x,y) < \frac{1-\alpha_n(x,y)-\beta_n(x,y)}{2\beta_n(x,y)}\right) \wedge (\beta_n(x,y) \geq 1)\right] \qquad (2.17)$$

However, $\beta_n > 1$ with $\xi_n > 0$ has to be excluded because of the unboundedness or non-negativity of the second right-hand-side term of (2.15).

B) $d(T^n x, T^n y) \geq d(x,y)$ for some $x, y \in X$, $n \in N$. Then, one gets from (2.10):

$$d(T^n x, T^n y) \leq (\alpha_n + \beta_n) d^2(x,y) + \beta_n d^2(T^n x, T^n y) + 2\mu_n \beta_n d^2(T^n x, T^n y) + \xi_n$$

$$\Rightarrow d(T^n x, T^n y) \leq k_{bn} d^2(x,y) + \frac{\xi_n}{1-\beta_n(1+2\mu_n)} \qquad (2.18)$$

where

$$k_{bn} = k_{bn}(x,y) = \frac{\alpha_n + \beta_n}{1-\beta_n(1+2\mu_n)} \to 1 \text{ as } n \to \infty \qquad (2.19)$$

which holds from (2.12) and $k_{bn} \geq 1$ if $\limsup_{n \to \infty} [\beta_n(x,y) \max(1, 1+2\mu_n(x,y))] < 1$, and

$$\mu_n(x,y) \in \left[\frac{1-\alpha_n(x,y)-2\beta_n(x,y)}{2\beta_n(x,y)}, \frac{1-\beta_n(x,y)}{2\beta_n(x,y)}\right) \qquad (2.20)$$

Thus, (2.15)-(2.16), with the second option in the logic disjunction being true if and only if $\xi_n = 0$ together with (2.18)-(2.20) are equivalent to (2.12)-(2.13) by taking $k_n = k_n(x,y)$ to be either $k_{an}$ or $k_{bn}$ for each $n \in N$. It then follows that $\exists \limsup_{n \to \infty} (d(T^n x, T^n y) - d(x,y)) \leq 0$; $\forall x, y \in X$ from (2.15)-(2.19) since $0 \leq k_n = k_n(x,y) \leq 1$ and $k_n(x,y) \to 1$; $\forall x, y \in X$ as $n \to \infty$. Thus, $T: X \to X$ is asymptotically nonexpansive. Hence, Property (i). Property (ii) is proven as follows. Consider the metric-induced norm $\|\ \|$, equivalent to the translation-invariant homogeneous metric $d: X \times X \to R_{0+}$. Such a norm exists since the metric is homogeneous and translation-invariant so that norm and metric are formally. Rename $A_0 \equiv A$ and define a sequence of subsets $A_j := \{T^j x : x \in A_0\}$ of $X$. From Property (i), $\{d(T^n x, T^n y)\}$ is bounded; $\forall x, y \in X$ if $d(x,y)$ is finite, since it is bounded for any finite $n \in N$ and, furthermore, it has a finite limit as $n \to \infty$. Thus, all the collections of subsets $\cup_{i=1}^k A_i$; $\forall k \in N$ are bounded since $A_0$ is bounded. Define the set $C = C(A_0) := cl\left[convex\left(\cup_{i=1}^\infty A_k\right)\right]$ which is nonempty bounded, closed and convex by construction. Since $(X,d)$ is complete, $(X, \|\ \|) \equiv (X,d)$ is a uniformly convex Banach space and $T: C \to C$ is asymptotically nonexpansive from Property (i) then, it has a fixed point $z = Tz \in C$, [1], [5]. Since the restricted self-mapping $T: C \to C$ is also continuous, one gets from Property (i):

$$\exists \lim_{n \to \infty} d(T^n x, T^n z) = \lim_{n \to \infty} d(T^n x, z) = d\left(\lim_{n \to \infty} T^n x, z\right) \leq d(x,z) < \infty \ ; \ \forall x \in A \qquad (2.21)$$

Then, any sequence $\{T^n x\}$ is convergent (otherwise, the above limit would not exist contradicting Property (i)), and then bounded, in $C$, $\forall x \in A$. This also implies $d(T^n x, T^{n+m} x)$ is bounded; $\forall x \in A$, $\forall n, m \in N$ and $\lim_{n \to \infty} d(T^n x, T^{n+m} x) = 0$; $\forall x \in A$, $\forall m \in N$. This implies also $T^n x \to z_x(x)$ as $n \to \infty$; $\forall x \in A$ such that $z_x(x) = Tz_x$; $\forall x \in A$ which is then a fixed point of $T: C \to C$ (otherwise, the above property $\lim_{n \to \infty} d(T^n x, T^{n+m} x) = 0$; $\forall x \in A$, $\forall m \in N$ would be contradicted). Hence, Property (ii). □

**Theorem 2.2**. Let $(X, d)$ be a metric space and consider the self-mapping $T: X \to X$. Assume that the constraint below holds:

$$d^2(T^n x, T^n y) \le \alpha_n(x, y) d^2(x, y) + \beta_n(x, y)(d^2(x, y) + d^2(T^n x, T^n y))$$
$$+ 2\mu_n(x, y)\beta_n(x, y) d^2(T^n x, T^n y) + \xi_n(x, y); \forall x, y \in X, \forall n \in N \quad (2.22)$$

with

$$\xi_n = \xi_n(x, y)$$
$$:= \max\left(0, (1 - \beta_n(x, y))d^2(T^n x, T^n y) - (\alpha_n(x, y) + \beta_n(x, y))d^2(x, y) - 2\mu_n(x, y)\beta_n(x, y)d^2(T^n x, T^n y)\right)$$
$$\to 0 \quad ; \forall x, y \in X \text{ as } n \to \infty \quad (2.23)$$

for some parameterizing real sequences $\alpha_n = \alpha_n(x, y)$, $\beta_n = \beta_n(x, y)$ and $\mu_n = \mu_n(x, y)$ satisfying for any $n \in N$:

$$\{\alpha_n(x, y)\} \subset [0, \infty), \{\mu_n(x, y)\} \subset \left[-1, \frac{1 - \beta_n(x, y)}{2\beta_n(x, y)}\right), \{\beta_n(x, y)\} \subset [0,1] ; \forall x, y \in X, \forall n \in N \quad (2.24)$$

Then the following properties hold:

(i) $\exists \lim_{n \to \infty} d(T^n x, T^n y) \le d(x, y)$ so that $T: X \to X$ is asymptotically nonexpansive, and then $\exists \lim_{n \to \infty} d(T^n x, T^n y) \le d(x, y)$; $\forall x, y \in X$, if

$$\frac{\alpha_n(x, y) + \beta_n(x, y)}{1 - \beta_n(x, y)(1 + 2\mu_n(x, y))} \ge 1 \Leftrightarrow \mu_n(x, y) \in \left[\frac{1 - \alpha_n(x, y) - 2\beta_n(x, y)}{2\beta_n(x, y)}, \frac{1 - \beta_n(x, y)}{2\beta_n(x, y)}\right); \forall x, y \in X, \forall n \in N$$

(2.25)

and the following limit exists:
$$\alpha_n(x, y) + 2\beta_n(x, y)(1 + \mu_n(x, y)) \to 1; \forall x, y \in X \text{ as } n \to \infty \quad (2.26)$$

(ii) Property (ii) of Theorem 2.1 if $(X, d)$ is complete and $(X, \|\ \|) \equiv (X, d)$ is a uniformly convex Banach under the metric-induced norm $\|\ \|$.

**Definition 2.3**. Assume that $(X,d)$ is a complete metric space with $d: X \times X \to R_{0+}$ being a homogeneous translation-invariant metric. Thus, $T: A \to A$ is asymptotically $\beta$-strictly pseudocontractive in the intermediate sense if:

$$\limsup_{n \to \infty} \left( (1 - \beta_n(1 + 2\mu_n))d^2(T^n x, T^n y) - (\alpha_n + \beta_n)d^2(x, y) \right) \leq 0; \quad \forall x, y \in A \quad (2.27)$$

for $\beta_n = \beta \in [0,1)$; $\forall n \in N$ and some real sequences $\{\alpha_n\}$, $\{\mu_n\}$ being, in general, dependent on the initial points, i.e. $\alpha_n = \alpha_n(x,y)$, $\mu_n = \mu_n(x,y)$ and:

$$\{\mu_n\} \subset \left[-1, \frac{1-\beta}{2\beta}\right) \text{ and } \{\alpha_n\} \subset [1, \infty); \forall n \in N, \alpha_n \to 1 \text{ and } \mu_n \to -1 \text{ as } n \to \infty \quad ; \forall x, y \in A, \forall n \in N$$

(2.28)

□

**Definition 2.4**. $T: A \to A$ is asymptotically pseudocontractive in the intermediate sense if (2.30) holds with $\{\mu_n\} \subset \left[-1, \frac{1-\beta_n}{2\beta_n}\right)$, $\{\beta_n\} \subset [0,1]$, $\{\alpha_n\} \subset [1, \infty)$, $\alpha_n \to 1$, $\beta_n \to 1$, $\mu_n \to -1$ as $n \to \infty$ and the remaining conditions as in Definition 2.3 with $\alpha_n = \alpha_n(x,y)$, $\beta_n = \beta_n(x,y)$ and $\mu_n = \mu_n(x,y)$. □

**Definition 2.5**. $T: A \to A$ is asymptotically $\beta$-strictly contractive in the intermediate sense if $\alpha_n \in [0, \infty)$, $\beta_n = \beta \in [0,1)$, $\mu_n \in \left[-1, \frac{1-\beta}{2\beta}\right); \forall n \in N$, $\mu_n \to \mu \in \left[-1, \frac{1-\beta}{2\beta} \min\left(1, \frac{1}{\alpha+\beta}\right)\right)$,

$\alpha_n \to \alpha \in [0,1)$ as $n \to \infty$, in Definition 2.3 with $\alpha_n = \alpha_n(x,y)$, $\mu_n = \mu_n(x,y)$. □

**Definition 2.6**. $T: A \to A$ is asymptotically contractive in the intermediate sense if $\alpha_n \in [0, \infty)$, $\{\beta_n\} \subset [0,1)$, $\mu_n \in \left[-1, \frac{1-\beta_n}{2\beta_n}\right)$ ; $\forall n \in N$, $\mu_n \to \mu \in \left[-1, -\frac{1+\alpha}{2}\right)$, $\alpha_n \to \alpha \in [0,1)$, and $\beta_n \to \beta = 1$ as $n \to \infty$ in Definition 2.3 with $\alpha_n = \alpha_n(x,y)$, $\beta_n = \beta_n(x,y)$ and $\mu_n = \mu_n(x,y)$. □

The next result relies on the asymptotically contractive and pseudocontractive mappings in the intermediate sense. Therefore, it is assumed that $\{\alpha_n(x,y)\} \subset [1, \infty)$.

**Theorem 2.8**. Let $(X, d)$ be a complete metric space endowed with a homogeneous translation-invariant metric $d: X \times X \to R_{0+}$ and consider the self-mapping $T: X \to X$. Assume that $(X, \|\ \|) \equiv (X, d)$ is a uniformly convex Banach space endowed with a metric-induced norm $\|\ \|$ from the metric $d: X \times X \to R_{0+}$. Assume that the asymptotically nonexpansive condition (2.22), subject to (2.23), holds for some parameterizing real sequences $\alpha_n = \alpha_n(x,y)$, $\beta_n = \beta_n(x,y)$ and $\mu_n = \mu_n(x,y)$ satisfying for any $n \in N$ :

$$\{\alpha_n(x,y)\} \subset [1, \infty), \{\mu_n(x,y)\} \subset \left[-1, \frac{1-\beta_n(x,y)}{2\beta_n(x,y)}\right), \{\beta_n(x,y)\} \subset [0, \beta) \subset [0,1] \quad (2.29)$$

$;\forall x, y \in X, \forall n \in N$. Then, $\exists \lim_{n \to \infty} d(T^n x, T^n y) \leq d(x, y)$ for any $x, y \in X$ satisfying the conditions:

$$\frac{\alpha_n(x, y) + \beta_n(x, y)}{1 - \beta_n(x, y)(1 + 2\mu_n(x, y)d(x, y))} \geq 1 \; ; \; \alpha_n(x, y) + 2\beta_n(x, y)(1 + \mu_n(x, y)) \to 1; \; \forall x, y \in X \text{ as } n \to \infty$$

(2.30)

Furthermore, the following properties hold:

**(i)** $T: C \to C$ is asymptotically $\beta$-strictly pseudocontractive in the intermediate sense for some nonempty, bounded, closed and convex set $C = C(A) \subset X$ and any given nonempty, bounded and closed subset $A \subset X$ of initial conditions if (2.29) hold with $0 \leq \beta_n = \beta < 1$, $\{\mu_n\} \subset \left[-1, \frac{1-\beta}{2\beta}\right)$, $\{\alpha_n\} \subset [1, \infty)$, $\alpha_n \to 1$ and $\mu_n \to -1$, as $n \to \infty$; $\forall x, y \in A$, $\forall n \in N$. Also, $T: C \to C$ has a fixed point, for any such a set $C$, if $T: X \to X$ is continuous.

**(ii)** $T: C \to C$ is asymptotically pseudocontractive in the intermediate sense for some nonempty, bounded, closed and convex set $C = C(A) \subset X$ and any given nonempty, bounded and closed subset $A \subset X$ of initial conditions if (2.29) hold with $\{\beta_n\} \subset [0, 1]$, $\{\mu_n\} \subset \left[-1, \frac{1-\beta_n}{2\beta_n}\right)$, $\{\alpha_n\} \subset [1, \infty)$, $\beta_n \to 1$, $\alpha_n \to 1$ and $\mu_n \to -1$, as $n \to \infty$; $\forall x, y \in A$, $\forall n \in N$. Also, $T: C \to C$ has a fixed point, for any such a set $C$, if $T: X \to X$ is continuous.

**(iii)** If (2.29) hold with $\alpha_n \in [0, \infty)$, $\beta_n = \beta \in [0, 1)$, $\mu_n \in \left[-1, \frac{1-\beta}{2\beta}\right)$, $\mu_n \to \mu \in \left[-1, \frac{1-\alpha-2\beta}{2\beta}\right)$; $\forall n \in N$ and $\alpha_n \to \alpha \in [0, 1)$ as $n \to \infty$ then $T: X \to X$ is asymptotically $\beta$-strictly contractive in the intermediate sense. Also, $T: X \to X$ has a unique fixed point.

**(iv)** If (2.29) hold with $\alpha_n \in [0, \infty)$, $\{\beta_n\} \subset [0, 1)$, $\mu_n \in \left[-1, \frac{1-\beta_n}{2\beta_n}\right)$, $\mu_n \to \mu \in \left[-1, -\frac{1+\alpha}{2}\right)$; $\forall n \in N$, $\beta_n \to 1$ and $\alpha_n \to \alpha \in [0, 1)$ as $n \to \infty$ then $T: X \to X$ is asymptotically strictly contractive in the intermediate sense. Also, $T: X \to X$ has a unique fixed point.

**3. Asymptotic contractions and pseudocontractions of cyclic self-mappings in the intermediate sense**

Let $A, B \subset X$ be nonempty subsets of $X$, $T: A \cup B \to A \cup B$ is cyclic self-mapping if $T(A) \subseteq B$ and $T(B) \subseteq A$. Assume that asymptotically nonexpansive condition (2.10), subject to (2.11), is modified as follows:

$$d^2(T^n x, T^n y) \leq \alpha_n(x, y) d^2(x, y) + \beta_n(x, y)(d^2(x, y) + d^2(T^n x, T^n y))$$

$$+2\mu_n(x,y)\beta_n(x,y)d(x,y)d(T^nx,T^ny)+\xi_n(x,y)+\gamma_n(x,y)D^2\;;\;\forall x\in A, y\in B,\forall n\in N \quad (3.1)$$

$$\xi_n = \xi_n(x,y)$$
$$:= max\Big(0,\,(1-\beta_n(x,y))d^2(T^nx,T^ny)-(\alpha_n(x,y)+\beta_n(x,y))d^2(x,y)-2\mu_n(x,y)\beta_n(x,y)d(x,y)d(T^nx,T^ny)\Big)$$
$$;\;\forall x\in A, y\in B,\forall n\in N \quad (3.2)$$

with $(\xi_n - \gamma_n(x,y)D^2)\to 0$; $\forall x,y \in X$ as $n\to\infty$, and that that asymptotically nonexpansive condition (2.22), subject to (2.23), is modified as follows:

$$d^2(T^nx,T^ny)\le \alpha_n(x,y)d^2(x,y)+\beta_n(x,y)\big(d^2(x,y)+d^2(T^nx,T^ny)\big)$$
$$+2\mu_n(x,y)\beta_n(x,y)d^2(T^nx,T^ny)+\xi_n(x,y)+\gamma_n(x,y)D^2\;;\;\forall x\in A, y\in B,\forall n\in N \quad (3.3)$$

$$\xi_n = \xi_n(x,y)$$
$$:= max\Big(0,\,(1-\beta_n(x,y))d^2(T^nx,T^ny)-(\alpha_n(x,y)+\beta_n(x,y))d^2(x,y)-2\mu_n(x,y)\beta_n(x,y)d(T^nx,T^ny)\Big)$$
$$;\;\forall x,y\in X \text{ as } n\to\infty \quad (3.4)$$

with $(\xi_n - \gamma_n(x,y)D^2)\to 0$; $\forall x,y\in X$ as $n\to\infty$, where $\{\gamma_n(x,y)\}\in[0,\infty)$ and $D = dist(A,B)\ge 0$. If $A\cap B \ne \emptyset$ then $D = 0$ and Theorems 2.1, 2.2 and 2.8 hold with the replacement $A \to A\cap B$. Then, if $A$ and $B$ are closed and convex then there is a unique fixed point of $T: A\cup B \to A\cup B$ in $A\cap B$. In the following, we consider the case that $A\cap B = \emptyset$ so that $D > 0$. The subsequent result based on Theorems 2.1, 2.2 and 2.8 holds:

**Theorem 3.1**. Let $(X,d)$ be a metric space and let $T: A\cup B \to A\cup B$ be a cyclic self-mapping, i.e. $T(A)\subseteq B$ and $T(B)\subseteq A$, where $A$ and $B$ are nonempty subsets of $X$. Define the sequence $\{k_n\}_{n\in N} \subset [0,\infty)$ of asymptotically nonexpansive iteration-dependent constants as follows:

$$k_n = k_n(x,y) := \begin{cases} \dfrac{\alpha_n + \beta_n(1+2\mu_n)}{1-\beta_n} \le 1 & \text{if } d(T^nx,T^ny)\le d(x,y) \quad (3.5a) \\[6pt] \dfrac{\alpha_n + \beta_n}{1-\beta_n(1+2\mu_n)} \ge 1 & \text{if } d(T^nx,T^ny)\ge d(x,y) \quad (3.5b) \end{cases}$$

$$;\forall (x,y)\in (A\times B)\cup(B\times A), \forall n\in N \quad \text{provided that } T: A\cup B\to A\cup B$$

satisfies the constraint (3.1), subject to (3.2), and

$$\Big[\big(d(T^nx,T^ny)\ge d(x,y)\wedge \beta_n = 1\big)\Rightarrow (\gamma_n = 0), \forall (x,y)\in (A\times B)\cup(B\times A),\forall n\in N\Big] \quad (3.6)$$

and

$$k_n = k_n(x,y) = \dfrac{\alpha_n + \beta_n}{1-\beta_n(1+2\mu_n)} \ge 1 \quad (3.7)$$

; $\forall n \in N$ for $x \in A (y \in B)$ and for $x \in B (y \in A)$ provided that $T: A \cup B \to A \cup B$ satisfies the constraint (3.3), subject to (3.4) provided that the parameterizing bounded real sequences $\{\alpha_n(x, y)\}, \{\beta_n(x, y)\}, \{\mu_n(x, y)\}$ and $\{\gamma_n(x, y)\}$ of general terms $\alpha_n = \alpha_n(x, y)$, $\beta_n = \beta_n(x, y)$ and $\mu_n = \mu_n(x, y)$ fulfil the following constraints:

$$\left[\left(\mu_n(x, y) \in \left[-\frac{\alpha_n(x, y) + \beta_n(x, y)}{2\beta_n(x, y)}, \frac{1 - \alpha_n(x, y) - 2\beta_n(x, y)}{2\beta_n(x, y) d(x, y)}\right]\right) \wedge (\beta_n(x, y) < 1)\right]$$

$$\vee \left[\left(\mu_n(x, y) < -\frac{\alpha_n(x, y) + \beta_n(x, y)}{2\beta_n(x, y)}\right) \wedge (\beta_n(x, y) > 1) \Leftrightarrow \xi_n(x, y) = 0\right]$$

$$\vee \left[\mu_n(x, y) \in \left[\frac{1 - \alpha_n(x, y) - 2\beta_n(x, y)}{2\beta_n(x, y)}, \frac{1 - \beta_n(x, y)}{2\beta_n(x, y)}\right)\right]; \forall (x, y) \in (A \times B) \cup (B \times A), \forall n \in N \quad (3.8)$$

, $\gamma_n = \gamma_n(x, y) \geq \max(0, 1 - k_n)$ and assuming that the following limits exist:

$$\left(\mu_n(x, y) - \frac{1 - \alpha_n(x, y) - 2\beta_n(x, y)}{2\beta_n(x, y)}\right) \to 0 \Leftrightarrow \alpha_n + 2\beta_n(1 + \mu_n) \to 1 ; \gamma_n(x, y) \to 0$$

$$; \forall (x, y) \in (A \times B) \cup (B \times A) \text{ as } n \to \infty \quad (3.9)$$

Then, the following properties hold:

(i) $T: A \cup B \to A \cup B$ satisfies (3.3), subject to (3.4)-(3.9); $\forall (x, y) \in (A \times B) \cup (B \times A)$). Then,

$$\exists \lim_{n \to \infty} d(T^n x, T^n y) \in [D, d(x, y)]; \forall (x, y) \in (A \times B) \cup (B \times A)$$

so that $T: A \cup B \to A \cup B$ is a cyclic asymptotically nonexpansive self-mapping. If $x \in A$ is a best proximity point of $A$ and $y \in B$ is a best proximity point of $B$ then $\lim_{n \to \infty} d(T^n x, T^n y) = D$ and $T^{2n} x \to z_x = z(x)$ and $T^{2n} y \to z_y = z(y)$ which are best proximity points of $A$ and $B$ (not being necessarily identical to $x$ and $y$), respectively if $T: A \cup B \to A \cup B$ is continuous.

(ii) Property (i) also holds if $T: A \cup B \to A \cup B$ satisfies (3.1), subject to (3.2), (3.7), (3.8)-(3.9) and (3.5b) provided that $d(T^n x, T^n y) \leq d(x, y); \forall (x, y) \in (A \times B) \cup (B \times A)$.

**Theorem 3.3.** Let $(X, d)$ be a metric space and let $T: A \cup B \to A \cup B$ be a cyclic self-mapping which satisfies the asymptotically nonexpansive constraint (3.1), subject to (3.2), where $A$ and $B$ are nonempty subsets of $X$. Let the sequence $\{k_n\}_{n \in N} \subset [0, 1)$ of asymptotically nonexpansive iteration-dependent constants be defined by a general term

$$k_n(x, y) = k_n := \frac{\alpha_n + \beta_n(1 + 2\mu_n)}{1 - \beta_n} \in [0, 1) \quad (3.10)$$

under the constraints

$$\gamma_n(x,y) = \gamma_n := \delta_j(1-k_n)(1-\beta_n) = o(1-\beta_n),\ \beta_n \leq 1 \Rightarrow \mu_n \leq -\frac{1+\alpha_n}{2},\ \forall n \in N \text{ and } \lim_{n\to\infty} k_n = 1 \quad (3.11)$$

Then, the subsequent properties hold:

**(i)** The following limits exist:

$$\lim_{n\to\infty} d(T^n x, T^n y) = D;\ \forall (x,y) \in (A\times B) \cup (B\times A);\ \lim_{n\to\infty} d(T^n x, T^{n+1} x) = D;\ \forall x \in A \cup B \quad (3.12)$$

**(ii)** Assume, furthermore, that $(X,d)$ is complete, $A$ and $B$ are closed and convex and $d: X\times X$ is translation-invariant and homogeneous and $(X,d) \equiv (X, \|\ \|)$ is uniformly convex where $\|\ \|$ is the metric-induced norm. Then,

$$\lim_{n\to\infty} d(T^{2n} x, T^{2n+2} x) = \lim_{n\to\infty} d(T^{2n+1} x, T^{2n+3} x) = 0;\ \forall x \in A \cup B \quad (3.13)$$

$\{T^{2n} x\} \to z \in A$, $\{T^{2n+1} x\} \to Tz \in B$; $\forall x \in A$, and $\{T^{2n} y\} \to Tz \in B$, $\{T^{2n+1} x\} \to z \in A$ ; $\forall x \in A$, $\forall y \in B$, where $z$ and $Tz$ are unique best proximity points of $T: A\cup B \to A\cup B$ in $A$ and $B$, respectively. If $A\cap B \neq \emptyset$ then $z = Tz$ is the unique fixed point of $T: A\cup B \to A\cup B$.

It is now assumed that the cyclic self-mapping $T: A\cup B \to A\cup B$ is asymptotically nonexpansive while not being strictly contractive for any finite number of iterations contrarily to the above result. The concepts of cyclic pseudocontractions and strict contraction in the intermediate sense play an important role in the obtained results.

**Theorem 3.5.** Let $(X, \|\ \|)$ be a uniformly convex Banach space endowed with a metric-induced norm $\|\ \|$ from a translation-invariant homogeneous metric $d: X\times X \to R_{0+}$, where $A$ and $B \subset X$ are nonempty, closed and convex subsets of $X$ and assume that $T: A\cup B \to A\cup B$ is cyclic self-mapping. Define the sequence $\{k_n\}_{n\in N} \subset [0,\infty)$ of asymptotically nonexpansive iteration-dependent constants as follows:

$$k_n = k_n(x,y) := \begin{cases} \dfrac{\alpha_n + \beta_n(1+2\mu_n)}{1-\beta_n} \leq 1 & \text{if } d(T^n x, T^n y) < d(x,y) \quad (3.14) \\ \dfrac{\alpha_n + \beta_n}{1-\beta_n(1+2\mu_n)} \geq 1 & \text{if } d(T^n x, T^n y) \geq d(x,y) \quad (3.15) \end{cases}$$

; $\forall (x,y) \in (A\times B)\cup(B\times A), \forall n \in N$ provided that $T: A\cup B \to A\cup B$ satisfies the constraint (3.1), subject to (3.2); and

$$k_n = k_n(x,y) = \frac{\alpha_n + \beta_n}{1-\beta_n(1+2\mu_n)} \geq 1 \quad (3.16)$$

; $\forall n \in N$ for $x \in A (y\in B)$ and for $x \in B(y\in A)$ provided that $T: A\cup B \to A\cup B$ satisfies the constraint (3.3), subject to (3.4), provided that the parameterizing bounded real sequences $\{\alpha_n(x,y)\}, \{\beta_n(x,y)\}, \{\mu_n(x,y)\}$ and $\{\gamma_n(x,y)\}$ of general terms $\alpha_n = \alpha_n(x,y)$, $\beta_n = \beta_n(x,y)$ and $\mu_n = \mu_n(x,y)$ fulfil the following constraints:

$$\left[\left(\mu_n(x,y)\in\left[-\frac{\alpha_n(x,y)+\beta_n(x,y)}{2\beta_n(x,y)},\frac{1-\alpha_n(x,y)-2\beta_n(x,y)}{2\beta_n(x,y)d(x,y)}\right]\right)\wedge(\beta_n(x,y)<1)\right]$$

$$\vee\left[\left(\mu_n(x,y)<-\frac{\alpha_n(x,y)+\beta_n(x,y)}{2\beta_n(x,y)}\right)\wedge(\beta_n(x,y)>1)\Leftrightarrow\xi_n(x,y)=0\right]$$

$$\vee\left[\mu_n(x,y)\in\left[\frac{1-\alpha_n(x,y)-2\beta_n(x,y)}{2\beta_n(x,y)},\frac{1-\beta_n(x,y)}{2\beta_n(x,y)}\right)\right];\forall(x,y)\in(A\times B)\cup(B\times A),\forall n\in N \quad (3.17)$$

, $\gamma_n=\gamma_n(x,y)\geq max(0,1-k_n)$ and assuming that the following limits exist:

$$\mu_n(x,y)\to\frac{1-\alpha_n(x,y)-2\beta_n(x,y)}{2\beta_n(x,y)}\Leftrightarrow\alpha_n+2\beta_n(1+\mu_n)\to 1\ ;\ \gamma_n(x,y)\to 0$$

$$;\forall(x,y)\in(A\cup B)\times(B\cup A)\text{ as }n\to\infty \quad (3.18)$$

Then, the following properties hold:

**(i)** If $T:A\cup B\to A\cup B$ satisfies (3.3), subject to (3.14)-(3.18); $\forall(x,y)\in(A\times B)\cup(B\times A))$. Then,

$$\exists\lim_{n\to\infty}d(T^n x,T^n y)\in[D,d(x,y)]\ ;\ \forall(x,y)\in(A\times B)\cup(B\times A) \quad (3.19)$$

so that $T:A\cup B\to A\cup B$ is asymptotically nonexpansive. If $x\in A$ is a best proximity point of $A$ and $y\in B$ is a best proximity point of $B$ then $\lim_{n\to\infty}d(T^n x,T^n y)=D$ and $T^{2n}x\to z_x=z(x)$ and $T^{2n}y\to z_y=z(y)$ which are best proximity points of $A$ and $B$ (not being necessarily identical to $x$ and $y$), respectively, if furthermore, $T:A\cup B\to A\cup B$ is continuous.

**(ii)** Property (i) also holds if $T:A\cup B\to A\cup B$ satisfies (3.1), subject to (3.2), (3.22), (3.17)-(3.18) and (3.5b), with $d(T^n x,T^n y)\leq d(x,y)$; $\forall(x,y)\in(A\cup B)\times(B\cup A)$.

**(iii)** Assume that $T:A\cup B\to A\cup B$ is asymptotically $\beta$-strictly pseudocontractive in the intermediate sense so that (3.21) holds with $0\leq\beta_n=\beta<1$, $\{\mu_n\}\subset\left[-1,\frac{1-\beta}{2\beta}\right)$, $\{\alpha_n\}\subset[1,\infty)$, $\alpha_n\to 1$ and $\mu_n\to -1$, as $n\to\infty$ ; $\forall x,y\in A$, $\forall n\in N$. Then, $T:A\cup B\to A\cup B$ is asymptotically nonexpansive and Property (i) holds.

**(iv)** $T:A\cup B\to A\cup B$ is asymptotically pseudocontractive in the intermediate sense if (3.22) holds with $\{\beta_n\}\subset[0,1]$, $\{\mu_n\}\subset\left[-1,\frac{1-\beta_n}{2\beta_n}\right)$, $\{\alpha_n\}\subset[1,\infty)$, $\beta_n\to 1$, $\alpha_n\to 1$ and $\mu_n\to -1$, as $n\to\infty$ ; $\forall x,y\in A$, $\forall n\in N$. Then, $T:A\cup B\to A\cup B$ is asymptotically nonexpansive and Property (i) holds.

**(v)** If the conditions of Property (iv) are modified as $\alpha_n \in [0, \infty)$, $\mu_n \in \left[-1, \frac{1-\beta}{2\beta}\right)$, $\mu_n \to \mu \in \left[-1, \frac{1-\alpha-2\beta}{2\beta}\right)$; $\forall n \in N$, $\alpha_n \to \alpha \in [0,1)$ as $n \to \infty$ and $\beta_n = \beta \in [0,1)$ in (3.22) then $T: A \cup B \to A \cup B$ is asymptotically $\beta$ - strictly contractive in the intermediate sense. Also, $T: A \cup B \to A \cup B$ has a unique best proximity point $z$ in $A$ and unique best proximity point $Tz$ in $B$ to which the sequences $\{T^{2n}x\}$ and $\{T^{2n+1}x\}$ converge; $\forall x \in A$. If $x \in B$ then $T^{2n}x \to Tz$ and $T^{2n+1}x$ as $n \to \infty$.

**(vi)** If (3.4) is modified by $\alpha_n \in [0, \infty)$, $\{\beta_n\} \subset [0,1)$, $\mu_n \in \left[-1, \frac{1-\beta_n}{2\beta_n}\right)$, $\mu_n \to \mu \in \left[-1, -\frac{1+\alpha}{2}\right)$, ; $\forall n \in N$, $\alpha_n \to \alpha \in [0,1)$ and $\beta_n \to 1$ as $n \to \infty$ then $T: A \cup B \to A \cup B$ is asymptotically - strictly contractive in the intermediate sense. Also, $T: A \cup B \to A \cup B$ has a unique best proximity point in $A$ and a unique best proximity point in $B$ to which the sequences $\{T^{2n}x\}$ and $\{T^{2n+1}x\}$ converge as in Property (v).

**ACKNOWLEDGMENTS**

The author is grateful to the Spanish Ministry of Education for its partial support of this work through Grant DPI2009-07197. He is also grateful to the Basque Government for its support through Grant IT378-10.